\documentclass[11pt,letterpaper]{amsart}

\usepackage[margin=1.1in]{geometry}
\usepackage[colorlinks=true, linkcolor=blue, citecolor=blue]{hyperref}
\usepackage{bookmark}
\usepackage{amsmath}
\usepackage{amscd,latexsym}
     \usepackage{amssymb,amsfonts,bm}
\usepackage[all,arc]{xy}
\usepackage{enumerate}
\usepackage{mathrsfs}
\usepackage{tikz}
\usepackage{tikz-cd}
\usepackage{amsmath}
\usepackage{mathdots}
\usepackage{amsrefs}
\usepackage{graphicx}
\usepackage{mathtools}
\usepackage{leftidx}
\usepackage{tensor}
\usepackage{bbm}
\usepackage{dsfont}
\usepackage{enumitem,xcolor}
\usepackage{bbm}
\usepackage{multicol}

\theoremstyle{plane}
\newtheorem{theorem}{Theorem}[section]
\newtheorem{lemma}[theorem]{Lemma}
 \newtheorem{corollary}[theorem]{Corollary}
      \newtheorem{proposition}[theorem]{Proposition}

        \newtheorem{theorem*}{Theorem}

    \newcounter{relctr} %% <- counter for relations
\everydisplay\expandafter{\the\everydisplay\setcounter{relctr}{0}} %% <- reset every eq
\AtBeginDocument{} %% <- store original definition

\theoremstyle{definition}
\newtheorem*{definition*}{Definition}

 \newtheorem{example}[theorem]{Example}

\theoremstyle{remark}
\newtheorem{remark}[theorem]{Remark}
\newtheorem*{acknowledgments}{Acknowledgments}

\numberwithin{equation}{section}

\tikzset{
op/.style = {box/.style = {draw, rounded corners, fill=blue!30},
             level distance=9mm, sibling distance=9mm,
             baseline = (current bounding box.center)}
        }

\usepackage{pgfplots}

\makeatletter
\@namedef{subjclassname@2020}{%
  \textup{2020} Mathematics Subject Classification}
\makeatother

\begin{document}
\title{Covexillary Schubert varieties and Kazhdan-Lusztig Polynomials}
\author{Minyoung Jeon}
\address{Department of mathematics, University of Georgia, Athens, GA 30602, USA}
\email{minyoung.jeon@uga.edu}
 %14B05, 05E14
%\date{\today}
%
\begin{abstract}
We establish combinatorial and inductive formulas for Kazhdan-Lusztig polynomials associated to covexillary elements in classical types, extending results of Boe, Lascoux-Sch\"{u}tzenberger, Sankaran-Vanchinathan, and Zelevinsky for Grassmannians of classical types. The proof uses intersection cohomology theory and the isomorphism of Kazhdan-Lusztig varieties from Anderson-Ikeda-Jeon-Kawago. 
\end{abstract}
\maketitle

\section{Introduction}

 The celebrated Kazhdan-Lusztig polynomials were introduced by Kazhdan and Lusztig in their study of representations of Hecke algebras \cite{KL0}. Later the Kazhdan-Lusztig polynomials $P_{v, w}(q)$ were realized as Hilbert-Poincar\'{e} polynomials for the stalk intersection cohomology sheaf of the Schubert variety $X_{w}$ at a torus-fixed point $p_{v}\in X_{ w}$. Here $w$ and $v$ are elements in the Weyl group $W$ of a semisimple, simply connected algebraic group $G$ over an algebraically closed field \cite{KL}. This realization is related to the singularities of Schubert varieties, since the Kazhdan-Lusztig polynomial $P_{v, w}(q)$ for $w, v\in W$ equals $1$ if and only if the torus-fixed point $p_{ v}$ is a rationally smooth point in the Schubert variety $X_{w}$. Besides Kazhdan-Lusztig polynomials, singularities of Schubert varieties have been actively studied for decades, see, e.g., \cite{BL}. Further it has been a long standing open problem to describe the coefficients of $P_{ v,w}(q)$ combinatorially. The study of Kazhdan-Lusztig polynomials is also significant in representation theory of Lie algebra $\mathfrak{g}$ of G, see \cite[pg. 76]{BL}. For more recent developments, and further references, see \cite{WY23}.
 
 Our main results give combinatorial rules for the coefficients of the Kazhdan-Lusztig polynomials associated to {\it covexillary elements}. Specifically, in this paper, we consider classical algebraic groups $G$ of rank $n$ over $\mathbb{C}$: $G=SL(n)$ for type $A$, $G=SO(2n+1)$ for type $B$, $G=Sp(2n)$ for type $C$ and $G=SO(2n)$ for type $D$.

In type $A$, $W$ may be identified with the symmetric group $S_n$ of permutations of $\{1,2,\ldots,n\}$. The {\it vexillary permutations} $w\in S_n$ are those that avoid the pattern $2143$; in other words, there is no sequence of four numbers $a<b<c<d$ with $w(c) < w(d) < w(a) < w(b)$. Billey and Lam \cite{BilleyLam} extended the notion of vexillary elements $w\in W$ from type $A$ to Weyl groups for the other classical types. Anderson and Fulton \cite{AF20} reformulated the definition of these vexillary elements with a {``triple"} in the geometric context of degeneracy loci. In this paper, we use the definition and properties of a vexillary (signed) permutation from the point of view in \cite{AF20}. 

Recall, each Weyl group $W$ is generated by a set of simple reflections $S=\{s_i:i=1,2,\ldots,r\}$ where $r$ is the rank of the associated Lie group $G$. The length $\ell(w)$ of $w\in W$ is the minimal length of any factorization of $w=s_{i_1}s_{i_2}\cdots s_{i_N}$ into simple reflections. Geometrically, $\ell(w)$ is the dimension of $X_w$ as an algebraic variety. Now, each $W$ has a unique longest length element, denoted $w_\circ$.

For a vexillary element $w$, we call $w_\circ w$ a {\emph{covexillary}} element. In general, the case of covexillary elements includes the case of {\it cograssmannian} elements, which are those with at most one right ascent, i.e., $\ell(ws)>\ell(w)$ for at most one simple reflection $s\in W$. The terminology ``cograssmannian" is borrowed from \cite{JW}.

 \begin{theorem}[Theorem \ref{M2}]\label{introtheorem}
Let ${w}$ be a covexillary element in $W$, and $v$ be any element that satisfies $w\geq v$. Then there are some cograssmannian elements $w_\lambda, v_\mu$ in some Weyl group ${W}$, with $w_\lambda\geq v_\mu$ such that
$$P_{{v},w}(q)=P_{ {v}_\mu, w_\lambda}(q).$$
\end{theorem}

In Section \ref{gss} we will explicitly describe $w_\lambda,w_\mu$ along with the Weyl group $W$. Our proof of Theorem \ref{introtheorem} uses the interpretation of the Kazhdan-Lusztig polynomials as Hilbert-Poincar\'{e} polynomials for intersection cohomology of Schubert varieties associated to covexillary Schubert variety $X_{ w}$ at a torus-fixed point $p_{ v}\in X_{ w}$. We use the local isomorphism between vexillary Schubert varieties and Grassmannian Schubert varieties, recently established in the author's joint work \cite[Theorem 7.1]{AIJK}. This local isomorphism \eqref{Equiv} induces a ring isomorphism of the respective intersection cohomologies (Theorem \ref{MainTh}), allowing us to extend results \cite{B88,LS,SV,Z} for cograssmannians Schubert varieties to covexillary cases. 

More precisely, we deduce combinatorial (Corollary \ref{MainI}) and inductive (Corollary \ref{MainII}) formulas computing the Kazhdan-Lusztig polynomial $P_{v,  w}$ for a covexillary element $w\in W$ and any $v$ with $ v\leq w$.

Our generalization of (co)grassmannian cases to (co)vexillary cases, outside the type $A$ result by Lascoux, is new. The combinatorial formulas (Corollary \ref{MainI}) are generalizations of the result by Lascoux and Sch\"{u}tzenberger for type $A$ and Boe for the other classical types, providing a new proof of Lascoux \cite{L} and an analogous result by Li and Yong \cite{LY} for type $A$ (co)vexillary Schubert varieties.

  Lascoux and Sch\"{u}tzenberger \cite{LS} developed an algorithmic way of computing the Kazhdan-Lusztig polynomials for a cograssmannian element ${w_\circ w_\nu}$ in terms of weighted trees for type $A$. Boe \cite{B88} has carried over Lascoux and Sch\"{u}tzenberger's formulas for the classical Lie types. Lascoux \cite{L} generalized the algorithm of \cite{LS} to the type A covexillary case. Using Lascoux's rule \cite{L}, Li and Yong \cite{LY} formulated another combinatorial rule for Kazhdan-Lusztig polynomials in type A, in terms of drift configurations. Our results are independent of \cite{L}.

Stimulated by Lascoux and Sch\"{u}tzenberger's result \cite{LS}, Zelevinsky \cite{Z} gave a geometric explanation and provided inductive formulas by constructing small resolutions for type $A$ Schubert varieties $X_{w_\lambda}$ associated to (co)grassmannian elements $w_\lambda,w_\mu$ \cite{Z}. Sankaran and Vanchinathan \cite{SV} obtained similar results for types $B$, $C$ and $D$ under certain conditions. 
Our inductive formula (Corollary \ref{MainII}) generalizes results of Zelevinsky for type $A$ and of Sankaran and Vanchinathan for types $B$, $C$ and $D$.

 The paper is organized as follows. Section \ref{secrev} defines the Kazhdan-Lusztig polynomials in terms of the stalk intersection cohomology sheaf of Schubert varieties. In Section \ref{sec4.1}, we discuss covexillary Schubert varieties and related notions that will be used for the formulas in Corollary \ref{MainI} and Corollary \ref{MainII}. In Section \ref{sec44} we provide the proof of Theorem \ref{introtheorem} as well as Corollaries \ref{MainI} and \ref{MainII}.

\section{Schubert varieties and Kazhdan-Lusztig polynomials}\label{secrev}
We provide notions and definitions for the Kazhdan-Lusztig polynomials in terms of the sheaf-theoretical {\it intersection cohomology} theory of Deligne-Goresky-MacPherson \citelist{\cite{GM} \cite{M}}, just as in \cite{KL}.

Let $G:=G(n)$ be a classical algebraic group of rank $n$ over $\mathbb{C}$. Thus, $G=SL(n)$ for type $A$, $G=SO(2n+1)$ for type $B$, $G=Sp(2n)$ for type $C$ and $G=SO(2n)$ for type $D$. We use the standard terminology that can be found in \cite{Brenti}. We fix a maximal torus $T$ of $G$. Let $B$ be a Borel subgroup of upper triangular matrices such that $T\subset B\subset G$, $B^{-}$ an opposite Borel subgroup with $B\cap B^{-}=T$ and $W_n:=N_T/T$ the Weyl group of $G$ for the normalizer $N_T$ of $T$. (We hope the standard use of $B$ for Borel subgoup and our use of $B$ for the type $B$ will not cause confusion. The distinction should be clear from context.)

\par
In fact, the Weyl group $W:=W_n$ of type $A$ is isomorphic to the symmetric group $S_n$. To be specific, $S_n$ can be identified with permutations of $1,\ldots,n$ for type $A$. The Weyl group $W$ of type $B$ and $C$ are regarded as a group of signed permutations $S_n \ltimes \{\pm 1\}^n$. Specifically, the Weyl group $W$ of type $B$ can be identified with permutations of 
\[
\overline{n},\ldots,\overline{1},0,1,\ldots, n.
\] Here the barred numbers indicate the negation, i.e., $\overline{a}=-a$. For type $C$, it is identified with permutations of $\overline{n},\ldots,\overline{1},1,\ldots, n$. As for type $D$, we take $W$ to be the subgroup of the signed permutations, a set of permutations of $\overline{n},\ldots,\overline{1},1,\ldots, n$ with even number of barred numbers. In this paper, we will use one-line notation to represent an element $w$ of the Weyl group $W$: we write $w=w(1)\;w(2)\;\cdots w(n)\in S_n$ in one-line notation for type $A$. Following \cite{AF20,AIJK}, the signed permutation $w$ of types $B$, $C$ and $D$ only lists the values on positive integers $w=w(1)\;w(2)\;\cdots w(n)$ in one-line notation, since $w(\overline{i})=\overline{w(i)}$ for each $i$. Here we take the natural order 
\[
\overline{n}<\overline{n-1}<\cdots<n-1<n.
\]

Let $Fl=G/B$ denote a flag variety. The {\it Bruhat decomposition} is $Fl=\underset{w\in W}{\bigsqcup} X_w^\circ$, where the {\it Schubert cell} $X_w^\circ:=BwB/B$ is locally closed subset isomorphic to the affine space of dimension $\ell(w)$. The {\it Schubert variety} $X_w$ is the closure of $X_w^\circ$ such that $X_w=\underset{v\leq w}{\bigsqcup}X_v^\circ$ where $\leq$ denotes Bruhat order on $W$. One characterization of Bruhat order is that $v\leq w$ in $W_n$ if and only if  a $T$-fixed point $p_{v}$ is in $X_{w}$.

Let $Y$ be a pure $n$-dimensional algebraic variety, admitting a {\it stratification} \cite{J1}. A stratification of $Y$ is a filtration 
\[
Y_0\subset Y_1\subset \cdots\subset Y_n=Y
\]
of $Y$ by closed subvarieties such that $Y_i- Y_{i-1}$ is smooth of pure dimension $i$ along which $Y$ is ``equisingular". The notion of equisingularity is a technical notion we will not need to inspect in our proofs. By letting $U_i=Y-Y_{n-i}$, we have the filtration
\[
\emptyset\subset U_1\subset \cdots\subset U_n.
\]
Let $j_i:U_i\hookrightarrow U_{i+1}$ be the inclusion. Then the {\it intersection cohomology (Deligne) complex} $\mathcal{IC}({Y})$ can be determined by the canonical isomorphism
\begin{equation}\label{def:IC}
\mathcal{IC}(Y)\cong \tau_{\leq n-1}Rj_{n*}\cdots\tau_{\leq 0}Rj_{1*}\underline{\mathbb{Q}}_{U_1}
\end{equation}
where $\underline{\mathbb{Q}}_{U_1}$ is the constant sheaf of rational numbers in degree $0$ on $U_1$, $\tau_{\leq k}$ is the truncation functor that replaces the sheaf cohomology by zero in degree above $k$, and $Rj_{i*}$ is higher direct image for all $i$.

Given any stratified variety $Y$, we let $\mathcal{H}^*(\mathcal{IC}(Y))$ denote the sheaf cohomology of the intersection cohomology complex $\mathcal{IC}({Y})$. We write $\mathcal{H}_y(\mathcal{IC}(Y))$ for the stalk of $\mathcal{H}(\mathcal{IC}(Y))$ at a point $y\in Y$. We note that Schubert varieties have the stratifications as the Whitney stratification: $X_w$ of dimension $\ell(w)$ is the union of Schubert cells $X_v^\circ$ for all $v\leq w$. Each Schubert cell is smooth of pure dimension and its closure is a union of strata with partial order on $W$.

\label{KLdef}
For $v\in W$ with $v\leq w$ the Kazhdan-Lusztig polynomial $P_{ v,w}(q)$ can be obtained by
\begin{equation}\label{Poincare}
P_{v,w}(q)=\sum_{i\geq 0}\mathrm{dim}\;\mathcal{H}^{2i}_{p_{ v}}(\mathcal{IC}(X_{w}))q^i,
\end{equation}
as the Hilbert-Poincar\'{e} polynomial for the stalk intersection cohomology sheaf of $X_{ w}$ at the $T$-fixed point $p_{ v}\in X_{ w}$  \cite[Theorem 4.3]{KL}.

\section{Covexillary Schubert varieties in classical types}\label{sec4.1}
Throughout the rest of this manuscript, we express covexillary elements in terms of vexillary element associated to a triple $\tau$ as in \cite{AF20}. Let $w_\circ$ be the longest element in $W$, more precisely, it takes  $i$ to $n+1-i$ for type $A$ and $w(i)$ to $\overline{w(i)}$ for types $B$, $C$ and $D$.
 In this section, we define covexillary Schubert varieties $X_{w_\circ w}$ for a vexillary element $w$ in each type.

\subsection{ Type $A$}
We start with a Type $A$ triple $\tau$ to generate a vexillary element $w(\tau)$ as the first step. Given a triple, the vexillary element arises uniquely as in \cite{AF20}.

A {\it triple} $\tau=(\bold{k},\bold{p},\bold{q})$ of length $s$ denotes three sequences 
$
 \bold{k}=(0<k_1<\cdots<k_s),$ $\bold{p}=( 0<p_1\leq\cdots\leq p_s\leq n)$ and $\bold{q}=(0< q_1\leq\cdots\leq q_s< n)
$
of integers
with a condition 
\begin{equation*} \label{ineq}
k_{i+1}-k_{i}< (q_{i+1}-q_{i})-(p_{i+1}-p_{i})
\end{equation*}
for $1\leq i<s$. 

\begin{proposition}
Given a triple $\tau=(\bold{k},\bold{p},\bold{q})$, there is a unique vexillary permutation $w(\tau)\in W$. The vexillary element $w(\tau)$ is minimal in Bruhat order in $S_n$ such that
$\#\{a\leq p_i\;|\; w(a)>n-q_i\}=k_i$
for all $i$ and $1\leq a\leq n$. The permutations arising in this way are vexillary and all vexillary permutations are from a unique triple. 
\end{proposition}
\begin{proof}
After replacing $q_i$ by $n-q_i$, we can construct a permutation $w(\tau)$ from a triple as in \cite[\S 3.1.1]{AIJK}. Then the statements follows by \cite[p.13]{AF12}, \cite{AF20} and \cite[\S 3.1]{AIJK}.
 \end{proof}

 Let us describe the construction of $w(\tau)$ from a triple as follows. We start by putting $k_1$ integers $n-q_1+1,\cdots, n-q_1+k_1$ in increasing order to the positions from $ p_1-k_1+1$ to $p_1$. Then insert $k_2-k_1$ different integers, specifically, $n-q_2+1,\ldots,n-q_2+k_2-k_1$ in increasing order, as far right as possible but to the left of position $p_2+1$; i.e., place $n-q_2+k_2-k_1$ at position $p_2$, and then place $n-q_2+k_2-k_1-1$ at position $p_2-1$ if it is available. If the position $p_2-1$ is already occupied, then place $n-q_2+k_2-k_1-1$ at position $p_2-2$, and then continue this process until all $k_2-k_1$ different integers are placed. 
Continue in this way until $k_s$ integers have been inserted. Then enter unused integers in the rest of positions in increasing order.

For instance, for the triple $\tau=(2\;3\;5\;6,2\;4\;5\;6,2\;4\;7\;8)$ for $n=9$, we recover the permutation as follows. First we put $8,9$ in positions $1$ and $2$, and then put $6$ in position $4$. Then place $3,4$ at the positions $3$ and $5$, and then put $2$ at position $6$. Then fill in the unused integers $1,5,7$ in positions $7,8,9$ to get $w(\tau)=8\;9\;3\;6\;4\;2\;1\;5\;7$ in $S_9$.

For the next step, we define a covexillary Schubert variety associated to the covexillary element $w_\circ w(\tau)$. 
Let $V$ be a $n$-dimensional vector space over $\mathbb{C}$. The flag variety $SL(n)/B$ is the variety $Fl^A(V)$ of complete flags
\[
E_\bullet:\{\mathbf{0}\}\subset E_1\subset E_2\subset \cdots\subset E_n=V
\]
where the dimension of $E_i$ is $i$ for $1\leq i\leq n$.

Given a fixed partial flag 
$$ F_{q_1}\subseteq\cdots\subseteq F_{q_s}\subseteq V$$ of $V$
 with $\mathrm{dim}(F_{q_i})=q_i$, the {\it covexillary Schubert variety} $X_{w_\circ w(\tau)}$ is given by  
\begin{equation}\label{eqn:3.1}
X_{w_\circ w(\tau)}=\{E_\bullet\;|\;\mathrm{dim}(E_{p_i}\cap F_{q_i})\geq k_i\;\;\text{for} \;\; 1\leq i\leq s\}\subseteq Fl^A(V).
\end{equation}
Here $X_{w_\circ w(\tau)}$ is the Zariski closure of the $B$-orbit of the T-fixed point $p_{w_\circ w(\tau)}$.

\subsubsection{Partition and $\alpha\beta$-word of a covexillary permutation}\label{sec2.1.1}
 
 For a vexillary permutation $w$ with a triple $\tau$, there is an associated partition $\nu:=\nu(\tau)=(\nu_1\geq \nu_2\geq\cdots\geq \nu_{k_s})$ such that
\begin{equation}\label{nuA}
\nu_{k_i}=n-q_i-p_i+k_i
\end{equation}
and $\nu_k=\nu_{k_i}$ for $k_{i-1}<k\leq k_i$ with $k_0=0$ for $1\leq i\leq s$. The partition $\lambda$ of a covexillary permutation $w_\circ w(\tau)$ is obtained by a partition inside the $n\times n$ rectangle such that 
$\lambda_i=n-\nu_{i}$ for $1\leq i\leq n.$

We define a function 
\begin{equation}\label{b}
\mathfrak{H}: \lambda\mapsto
\begin{pmatrix}
a_1 & \cdots &a_d \\
b_0& \cdots &b_{d-1}\\
\end{pmatrix}_\lambda
\end{equation} where 
$a_i$ indicates the number of repeated $\lambda_i$ in $\lambda=\lambda_1^{a_1}\lambda_{2}^{a_2}\cdots \lambda_{d}^{a_d}$ and $b_0=\lambda_1$, $b_i=\lambda_{i+1}-\lambda_i$ for $1\leq i<d$. In this way, the partition $\lambda$ determines a $2\times d$ matrix.

For each covexillary element $w_\circ w(\tau)\in W_n$ with a triple $\tau$, we associate $w_\circ w(\tau)$ with a word in two two symbols $\alpha$ and $\beta$ by 
 \[
\lambda(w_\circ w(\tau))= \beta^{b_0}\alpha^{a_1}\cdots\beta^{b_{d-1}}\alpha^{a_d}.  \]

The following is an example of finding the partition of a covexillary permutation $w_\circ w(\tau)$.
\begin{example}
Let $n=8$. Let $w(\tau)=1\;4\;7\;5\;2\;3\;6\;8$ in $S_8$ be a vexillary permutation with a triple $\tau=(1\;3,3\;4,2\;5)$. The partition associated to $w(\tau)$ is $\nu=(4,2,2)$. And, the partition $\lambda$ of the covexillary permutation $w_\circ w(\tau)$ is $(4,6,6,8,8,8,8,8)$ from $\lambda_1=4, \lambda_2=6,$ and $\lambda_3=8$.

Furthermore, we have
$
\mathfrak{H}(\lambda)=\begin{pmatrix}
1 & 2 &5 \\
4& 2 &2\\
\end{pmatrix}_\lambda,
$
and $\lambda(w_\circ w(\tau))= \beta^{4}\alpha^{1}\beta^{2}\alpha^{2}\beta^{2}\alpha^{5}$.

 \end{example}

 \subsubsection{Partition and $\alpha\beta$-word of a pair}
 
Let $w(\tau)$ be a vexillary permutation in $S_n$ with a triple $\tau=(\bf{k},\bf{p},\bf{q})$. We may write $w$ for $w(\tau)$ for simplicity. Let $v$ be any element in $S_n$ such that $w\leq v$ in Bruhat order. 

We define a weakly increasing sequence $\bold{k}'= (0<k_1'\leq \cdots\leq k_s')$
defined by
\[
k_i'=\#\{a\leq p_i\;|\;v(a)>n-q_i\}.
\]
It is worthwhile to note that $\bold{p}=( 0<p_1\leq\cdots\leq p_s\leq n)$ and $\bold{q}=(0< q_1\leq\cdots\leq q_s< n)$ depend on $w(\tau)$ while $\bf{k}'$ relies on $v$. Since $v\geq w$ in Bruhat order, we have   
$
k_i'\geq k_i.
$
See \cite[\S 3.1.3]{AIJK} for this property.

A {\it weak triple} $\tau'(w,v)=(\bold{k}',\bold{p},\bold{q})$ consists of the three sequences $\bf{k}'$ from $v$ and $\bf{p},\bf{q}$ from $w(\tau)$, which satisfies 
$$k'_{i+1}-k'_{i}\leq (q_{i+1}-q_{i})-(p_{i+1}-p_{i})$$
for $1\leq i<s$. To construct a (strong) triple without the equality, we remove $k'_{i}, p_{i}, q_{i}$ if the equality holds and $k'_{i+1},p_{i+1},q_{i+1}$ if $k'_{i}=k'_{i+1}$ from $\tau'$ to have another triple $\tau''(w,v)=(\bold{k}'',\bold{p}'',\bold{q}'')$. In other words, 
$$ 
(\bold{k}'',\bold{p}'',\bold{q}''):=\begin{cases}
(\cdots\leq k_{i-1}'\leq k_{i+1}'\leq \cdots ,\cdots\leq p_{i-1}\leq p_{i+1}\leq\cdots, \cdots\leq q_{i-1}\leq q_{i+1}\leq \cdots)\\
\quad\quad\quad\quad\quad\quad\quad\quad\quad\quad\quad\quad\quad\quad\quad\quad\quad\quad\text{if}\;k'_{i+1}-k'_{i}= (q_{i+1}-q_{i})-(p_{i+1}-p_{i}),\\
(\cdots\leq k_{i}'< k_{i+2}'\leq \cdots ,\cdots\leq p_{i}\leq p_{i+2}\leq\cdots, \cdots\leq  q_{i}\leq q_{i+2}\leq \cdots)\\
\quad\quad\quad\quad\quad\quad\quad\quad\quad\quad\quad\quad\quad\quad\quad\quad\quad\quad \text{if}\; k_i'=k_{i+1}'.\\
\end{cases}
$$
We replace $\tau'(w,v)$ by $\tau''(w,v)$ and will not distinguish $\tau'(w,v)$ and $\tau''(w,v)$.

The partition $\mu$ of a pair $(w_\circ w, w_\circ v)$ with a weak triple $\tau':=\tau'(w,v)$ is obtained as the partition $\lambda$ of $w_\circ w$ with a triple $\tau(w)$: have a partition $\xi:=\xi(\tau')$ by setting $\xi_{k_i'}=n-q_i-p_i+k_i'$ and $\xi_k=\xi_{k_i'}$ for $k_{i-1}'<k\leq k_i'$. The partition $\mu$ of a pair $(w_\circ w,w_\circ v)$ is obtained by taking the partition inside $n\times n$ rectangle such that $\mu_i=n-\xi_i$ for $1\leq i\leq n$.

We have the $2\times d$ matrix corresponding to the partition $\mu:=\mu(w_\circ w,w_\circ v)$ with $\tau'$ as the image of the partition $\mu$ under the function 
\[
\mathfrak{H}:\mu\mapsto \begin{pmatrix}
e_1 & \cdots &e_d \\
f_0& \cdots &f_{d-1}\\
\end{pmatrix}_{\mu}
\]
defined in the same way as \eqref{b} with a partition $\mu$: $a_i$ indicates the number of repeated $\mu_i$ in $\mu=\mu_1^{e_1}\mu_{2}^{e_2}\cdots \mu_{d}^{e_d}$ and $f_0=\mu_1$, $f_i=\mu_{i+1}-\mu_i$ for $1\leq i<d$.
From $e_i$ and $f_i$, we have a $\alpha\beta$-word 
\[
\mu(w_\circ w, w_\circ v)= \beta^{f_0}\alpha^{e_1}\cdots\beta^{f_{d-1}}\alpha^{e_d}.
\]
associated to the pair $(w_\circ w,w_\circ v)$.

The following is an example of getting the partition $\mu$ from a pair of a covexillary permutation $w_\circ w$ and $w_\circ v$ for any $v\geq w$ in $S_n$ as well as its $\alpha\beta$-word.
\begin{example}
Let $n=8$. Let $w(\tau)=1\;4\;7\;5\;2\;3\;6\;8$ in $S_8$ be a vexillary permutation with a triple $\tau=(1\;3,3\;4,2\;5)$. Let $v=8\;7\;6\;5\;4\;3\;2\;1\geq w(\tau)$. Then we get a weak triple $\tau'=(2\;4,3\;4,2\;5)$ such that $\mu=(3,3,5,5,8,8,8,8)$. Thus, we have
\[
\mathfrak{H}(\mu)=\begin{pmatrix}
2& 2 &4 \\
3& 2 &3\\
\end{pmatrix}_{\mu}\quad\text{and}\quad \mu(w_\circ w,w_\circ v)=\beta^{3}\alpha^{2}\beta^{2}\alpha^{2}\beta^{3}\alpha^{4}
\]

\end{example}

 \subsection{Type $C$}
Let us move on to the type $C$ covexillary cases with similar steps as type $A$. We begin with a type $C$ triple $\tau$ to get a (signed) vexillary element $w:=w(\tau)$.

A {\it triple} $\tau=(\bf{k},\bf{p},\bf{q})$ of type $C$ is a collection of $\bold{k}=(0<k_1<\cdots<k_s)$, $\bold{p}=(0< p_1\leq\cdots\leq p_s\leq n)$ and $\bold{q}=(0<q_1\leq\cdots\leq q_s\leq n)$,
requiring 
\begin{equation}\label{eqn2.3}
k_{i+1}-k_{i}< (p_{i+1}-p_{i})+(q_{i+1}-q_{i})
\end{equation}
 for $1\leq i< s$.

 Given a triple $\tau$, one can determine a signed permutation $w=w(\tau)\in W_n$, in one-line notation, such that $w$ is minimal with $\#\{a\geq n+1-p_i\;|\;w(a)\leq \overline{n+1-q_i}\}=k_i$ for all $i$, see \cite[p.16]{AF12} and \cite{AF20}. Specifically, we begin with inserting a sequence of $k_1$ consecutive barred elements 
 \[
 \overline{n-q_1+k_1},\overline{n-q_1+k_1-1},\ldots,\overline{n-q_1+2},\overline{n-q_1+1}
 \]
  in positions $n-p_1+1$ to $n-p_1+k_1$. Then insert $k_2-k_1$ barred numbers, bars of the smallest remaining integers $a_{k_2-k_1}>\cdots>a_1\geq  n-q_2+1$, placed in order (bars of decreasing sequences, $\overline{a_{k_2-k_1}},\ldots,\overline{a_1}$) at or to the right of position $n-p_2+1$ but as close as possible to that position. Repeat this until one places $k_s$ barred integers. Then place the remaining positive integers in increasing order in the unfilled positions.

\begin{example}
Let $n=8$ and $\tau=(2\;4\;5\;6,3\;4\;6\;8,3\;6\;6\;7)$. Then we put $\overline{7}$ and $\overline{6}$ in positions $6,7$. Then place $\overline{4},\overline{3}$ in positions $5$ and $8$. Then put $\overline{5}$ at position $3$. Then put $\overline{2}$ at position $1$. Then we place remaining positive integers $1,8$ in the unfilled positions $2$ and $4$. In this way, we get a permutation $w(\tau)=\overline{2}\;1\;\overline{5}\;8\;\overline{4}\;\overline{7}\;\overline{6}\;\overline{3}$ in $W_8$.
\end{example}

A signed permutation $w$ in the Weyl group $W$ is {\it vexillary} if $w=w(\tau)$ for some triple $\tau=(\bf{k},\bf{p},\bf{q})$, as in type A. In particular, the ``empty" triple is a triple so that the identity element is vexillary {\cite[Definition 2.1]{AF20}}.

From a vexillary signed permutation, one can reconstruct the triple as follows. $n+1-p_1$ is the position after the last descent. We start at that position and find a sequence of $k_1$ consecutive barred integers ending in $\overline{n+1-q_1}$. Then continue this process on the remaining, but considering elements as being in sequence if the only missing integers have already been taken cared of, see \cite[p.16]{AF12} and \cite{AF20} for more details.

We remark that there is a different definition of the vexillary elements from \cite[Theorem 2]{Tam}. But our definition of vexillary permutations is derived from the ones presented in \cite{AF12,AF20,AIJK}, and we do not use properties deduced from \cite{Tam}.

Next, we define a covexillary Schubert variety associated to a covexillary element $w_\circ w$. Let $V$ be a vector space of dimension $2n$ over $\mathbb{C}$ equipped with a nondegenerate skew-symmetric form. A vector space $E$ is {\it isotropic} if $E \subset E^\perp$ with respect to the given form. An {\it isotropic (complete) flag} is a filtration 
\[
E_\bullet:\{\mathbf{0}\}\subset E_1\subset E_2\subset \cdots\subset E_n\subset V\cong \mathbb{C}^{2n},
\]
of linear subspaces
where $E_i$ are isotropic with $\mathrm{dim}(E_i)=i$. The variety $Fl^C(V)=Sp(2n)/B_n$ is identified with the set of isotropic flags. 

We fix an isotropic partial flag 
$$ F_{q_1}\subseteq\cdots\subseteq F_{q_s}\subseteq V$$
 where the subscripts indicate their dimension. A {\it covexillary Schubert variety} $X_{w_\circ w}$ associated to a triple $\tau=(\bf{k},\bf{p},\bf{q})$ is given by
$$X_{w_\circ w}=\{E_\bullet\;|\;\mathrm{dim}(E_{p_i}\cap F_{q_i})\geq k_i\;\;\text{for} \;\; 1\leq i\leq s\}\subset Fl^C(V)$$
 where $\mathrm{dim}(E_{p_i})=p_i$, as described in type $A$.

\subsubsection{Partition and $\alpha\beta$-word of a covexillary permutation}
For each covexillary signed permutation $w_\circ w(\tau)$ with a triple $\tau$, there is the largest strict partition $\lambda(\tau)$ such that 
\[
\lambda_{k_i}=p_i+q_i\;\;\text{for}\;1\leq i\leq s.
\]
That is, $\lambda(\tau)=(0<\lambda_1<\cdots<\lambda_{k_s})$ has parts $\lambda_k=p_i+q_i-k_i+k$ whenever $k_{i-1}\leq k<k_i$.

We recall the notation of \cite{SV}. Given a partition $\lambda$, a maximal subsequence of consecutive integers in this partition $\lambda$ is referred to as a {\it block} of $\lambda$. The partition $\lambda$ is the concatenation of its blocks. Let $a_i$ be the length of the $i$th block. Then $\lambda$ determines a $2\times s$ matrix as the image of the following function
 \begin{equation}\label{eqn2.4}
\mathfrak{H}:\lambda\mapsto \begin{pmatrix}
a_1 & \cdots &a_s\\
b_0 & \cdots & b_{s-1}\\
\end{pmatrix}_\lambda
\end{equation}
where the image is obtained by
\begin{enumerate}
\item $b_0=\lambda_{a_1}-a_1$,
\item and $b_i=\lambda_{k_{i+1}}-(a_1+\cdots+a_{i+1})-(b_0+\cdots+b_{i-1})$ for $1\leq i\leq s-1$, $k_i=a_1+\cdots+a_i$.
 \end{enumerate}

 To each covexillary element $w_\circ w(\tau)\in W_n$ with a triple $\tau$, we associate $w_\circ w(\tau)$ with a word in two two symbols $\alpha$ and $\beta$ by 
 \[
\lambda(w_\circ w(\tau))= \beta^{b_0}\alpha^{a_1}\cdots\beta^{b_{s-1}}\alpha^{a_s}\quad\text{or}\quad \beta^{b_0}\alpha^{a_1}\cdots\beta^{b_{s-1}}\alpha^{a_s}\beta^{a_{s+1}}\]
  for $\sum_i a_i+\sum_j b_j=2n$.

 \begin{example}\label{Aex1}
Let $w(\tau)=5\;\overline{4}\;\overline{3}\;6\;\overline{1}\;\overline{2}\;7$ be a vexillary permutation in $W_7$ with a triple $\tau=(1\;2\;4,2\;3\;6,6\;7\;7)$.

The partition $\lambda$ associated to the covexillary (signed) permutation $w_\circ w(\tau)$ is $(8,10,12,13)$ via $\lambda_1=8, \lambda_2=10,$ and $\lambda_4=13$.
Moreover, we have
$$\mathfrak{H}(\lambda)=\begin{pmatrix}
1 &1& 2\\
7 & 1 &1\\
\end{pmatrix}_\lambda,\quad\text{and}\quad \lambda(w_\circ w(\tau))= \beta^{7}\alpha^{1}\beta^{1}\alpha^{1}\beta^{1}\alpha^{2}\beta.$$
\end{example}

 \subsubsection{Partition and $\alpha\beta$-word of a pair}
Given a vexillary signed permutation $w:=w(\tau)$ with a triple $\tau$, we take any $v$ such that $v\geq w$ in $W_n$. As in type $A$, there is a {\it weak triple} $\tau'=(\bold{k}',\bold{p},\bold{q})$, obtained from $k'_i=\#\{a\geq n+1-p_i\;|\;v(a)\leq \overline{n+1-q_i}\}$ such that it satisfies
\begin{equation}\label{equalweek}
k'_{i+1}-k'_{i}\leq (p_{i}-p_{i+1})+(q_{i}-q_{i+1})
\end{equation}
for $i=1,\ldots, s-1$. We use the same convention as before. That is, we delete $k'_{i}, p_{i}, q_{i}$ if the equality holds and $k'_{i+1},p_{i+1},q_{i+1}$ if $k'_{i}=k'_{i+1}$ from $\tau'$. 

The partition $\mu$ of a pair $(w_\circ w,w_\circ v)$ with a weak triple $\tau'$ is given by 
$
\mu_{k_i'}=p_i+q_i
$
and $\mu_k=p_i+q_i-k_i'+k$ whenever $k_{i-1}'\leq k<k_i'$.

For a partition $\mu$, we have the $2\times s$ matrix 
$
\mathfrak{H}(\mu)=\begin{pmatrix}
e_1 & \cdots &e_s\\
f_0 & \cdots & f_{s-1}\\
\end{pmatrix}_{\mu}
$
 from a function $\mathfrak{H}$ defined as \eqref{eqn2.4}.

 For each pair $(w_\circ w,w_\circ v)$ with a weak triple $\tau'$, the $\alpha\beta$-word associated to the pair is given by
 \[
 \mu(w_\circ w,w_\circ v)=\beta^{f_0}\alpha^{e_1}\cdots\beta^{f_{s-1}}\alpha^{e_s}\quad\text{or}\quad \beta^{f_0}\alpha^{e_1}\cdots\beta^{f_{s-1}}\alpha^{e_s}\beta^{e_{s+1}}\]
  for $\sum_i e_i+\sum_j f_j=2n$.

\begin{example}
Let $w(\tau)=5\;\overline{4}\;\overline{3}\;6\;\overline{1}\;\overline{2}\;7$ be a vexillary permutation in $W_7$ with a triple $\tau=(1\;2\;4,2\;3\;6,6\;7\;7)$. Let $v=7\;\overline{6}\;\overline{5}\;\overline{4}\;\overline{3}\;\overline{2}\;\overline{1}\in W_7$. Then we have a weak triple $\tau'=(1\;3\;6,2\;3\;6,6\;7\;7)$ associated to the pair $(w_\circ w(\tau),w_\circ v)$.

The partition $\mu$ of the pair is given by $(8, 9,10,11,12,13)$ from $\mu_1=8,\mu_3=10,\mu_6=13$, and 
\[
\mathfrak{H}(\mu)=\begin{pmatrix}
1& 2& 3\\
7 &0 &0\\
\end{pmatrix}_{\mu}\quad\text{and}\quad\mu(w_\circ w,w_\circ v)=\beta^{7}\alpha^{1}\beta^{0}\alpha^{2}\beta^{0}\alpha^{3}\beta=\beta^{7}\alpha^{6}\beta.
\]

\end{example}

\subsection{Type $D$}

We follow similar procedures as before.

A type $D$ triple $\tau=(\bf{k},\bf{p},\bf{q})$ is $0<k_1<\cdots<k_s,$ $0\leq p_1\leq\cdots\leq p_s< n,$ and $0\leq q_1\leq\cdots\leq q_s< n$ with the same inequality condition \eqref{eqn2.3} in type $C.$ In particular $\tau$ is obtained by subtracting $1$ from each $p_i$ and $q_i$ of the sequences in type $C$. If $k_s$ is odd, the triple is replaced by one with $k_{s+1} = k_s + 1$ and
$p_{s+1} = q_{s+1} = n$.

A (signed) vexillary permutation $w:=w(\tau)$ in $W_n$ is characterized as minimal in Bruhat order such that $\#\{a\geq n-p_i\;|\;w(a)\leq \overline{n-q_i}\}=k_i$ for all $i$. In particular, one can construct a vexillary signed permutation $w(\tau)$ for type $D$ as follows. We get a triple $\tau^+$ of type $C$ uniquely by adding $1$ to each $p_i$ and each $q_i$. Then we construct $w(\tau)=w(\tau^+)$ with $\tau^+$ from type $C$. 
\begin{example}
For $\tau=(1\;2\;4,0\;2\;5,0\;2\;5)$ we have $w=\overline{2}\;\overline{1}\;3\;\overline{4}\;5\;\overline{6}$ in $W_6, n=6$.
\end{example}

A vector space $V$ has dimension $2n$ over $\mathbb{C}$ and is equipped with a symmetric bilinear form. An {\it isotropic flag} is a sequence of subspaces
\[
E_\bullet:\{\mathbf{0}\}\subset E_1\subset E_2\subset \cdots\subset E_{n-1}\subset E_n\subset V\cong \mathbb{C}^{2n}
\]
where all $E_i$ are isotropic with respect to the form, i.e., $E_i\subset E_i^\perp$, and has dimension of $i$. We note that $E_n$ is a maximal isotropic subspace of $V$. The variety $Fl^D(V)=SO(2n)/B$ consists of the set of all isotropic flags.

We take a partial flag of isotropic subspaces
$$F_{q_1}\subseteq\cdots\subseteq F_{q_s}\subseteq F_n\subseteq V$$
for $\mathrm{dim}(F_{q_i})=q_i$. In fact $F_n$ is the maximal isotropic subspace of $V$. The {\it covexillary Schubert variety} $X_{w_\circ w}$ is obtained by the closure of the Schubert cell
$$X_{w_\circ w}^\circ=\{E_\bullet\;|\;\mathrm{dim}(E_{p_i}\cap F_{q_i})= k_i\;\;\text{for} \;\; 1\leq i\leq s\}\subseteq Fl^D(V)$$
with $\mathrm{dim}(E_{p_i})=p_i$ for all $i$. The condition $p_s+q_s=2n$ can be contained in the definition of $X_{w_\circ w}$ so that $\mathrm{dim}(E_n\cap F_n)\equiv n$ $(\mathrm{mod}\;2)$ holds on $X_{w_\circ w}$.  

\subsubsection{Partition and $\alpha\beta$-word of a covexillary permutation}
Given a signed vexillary permutation $w:=w(\tau)$ with a triple $\tau=(\bf{k},\bf{p},\bf{q})$, the partition $\lambda$ of the covexillary (signed) permutation $w_\circ w$ is given by the largest strict partition 
such that 
\[
\lambda_{k_i}=p_i+q_i+1\;\;\text{for}\;1\leq i\leq s.
\]
and $\lambda_k=p_i+q_i-k_i+k$ whenever $k_{i-1}\leq k<k_i$. If $p_s=q_s=n$, then $\lambda_{k_s}=2n$.
Alternatively, we have a triple $\tau^+$ by adding $1$ to each $p_i$ and each $q_i$. Then the partition $\lambda^+=(\lambda_1+1,\ldots,\lambda_s+1)$ is the strict partition from $(\tau')^+$ for type $C$. 

The partition $\lambda$ determines a $2\times d$ matrix $\mathfrak{H}(\lambda)$, and have a $\alpha\beta$-word $\lambda(w_\circ w)$ associated to $w_\circ w$ as before in type $C$.

\begin{example}
Let $w(\tau)=\overline{2}\;\overline{1}\;3\;\overline{4}\;5\;\overline{6}$ be a (signed) vexillary permutation in $W_6, n=6$ with $\tau=(1\;2\;4,0\;2\;5,0\;2\;5)$. Then $\lambda=(1,5,10,11)$ from $\lambda_1=1$, $\lambda_2=5$, and $\lambda_4=11$. 

Also, we have
\[
\mathfrak{H}(\lambda)=\begin{pmatrix}
1& 1& 2\\
1 &2 &4\\
\end{pmatrix}_\lambda\quad{and}\quad\lambda(w)=\beta^{1}\alpha^{1}\beta^{2}\alpha^{1}\beta^{4}\alpha^{2}\beta
\]

\end{example}
 
 \subsubsection{Partition and $\alpha\beta$-word of a pair}
Given a tripe $\tau$ with a vexillary (signed) permutation $w:=w(\tau)$, and any $v\geq w$, the weak triple  $\tau'=(\bold{k}', \bold{p},\bold{q})$ associated to a pair $(w_\circ w,w_\circ v)$ is obtained by $k_i'=\#\{a\geq n-p_i\;|\;v(a)\leq \overline{n-q_i}\}$.

The partition $\mu=(\mu_1,\ldots,\mu_s)$ associated to the pair can be determined as before with the weak triple $\tau'$.  Alternatively, we may get the partition $\mu^+((\tau')^+)=(\mu_1+1,\ldots,\mu_s+1)$ for type $C$.
As in type $C$, we have $2\times s$ matrices $\mathfrak{H}(\mu)$ and $\alpha\beta$-word $\mu(w_\circ w,w_\circ v)$ with the weak triple $(\tau')^+$.

\begin{example}
Let $w(\tau)=\overline{2}\;\overline{1}\;3\;\overline{4}\;5\;\overline{6}$ be a (signed) vexillary permutation in $W_6$ with $\tau=(1\;2\;4,0\;2\;5,0\;2\;5)$. Let $v=\bar{3}\;\bar{2}\;\bar{1}\;\bar{5}\;\bar{4}\;\bar{6}$. The weak triple $\tau'$ associated to $(w_\circ w(\tau),w_\circ v)$ is given by $(1\;3\;6,0\;2\;5,0\;2\;5)$. 

The partition $\mu=(1,4,5,9,10,11)$ is obtained from $\mu_1=1,\mu_3=5$ and $\mu_6=11$. In addition, we have
\[
\mathfrak{H}(\mu)=\begin{pmatrix}
1& 2& 3\\
1 &1 &3\\
\end{pmatrix}_{\mu}\quad\text{and}\quad\mu(w_\circ w(\tau),w_\circ v)=\beta^{1}\alpha^{1}\beta^{1}\alpha^{2}\beta^{3}\alpha^{3}\beta.
\]
\end{example}

\subsection{Type $B$}
This section closely resembles type $C$.
Let $V$ be a $2n+1$ dimensional vector space equipped with a symmetric bilinear form on it. The isotropic flag variety $Fl^B(V)=SO(2n+1)/B$ is identified with the set of isotropic flags
\[
E_\bullet:\{\mathbf{0}\}\subset E_1\subset \cdots E_{n-1}\subset E_n\subset V\cong\mathbb{C}^{2n+1}
\]
where $E_i$ are isotropic with $\mathrm{dim}(E_i)=i$. Especially, $E_n$ is a maximal isotropic subspace of $V$.

We use the same triple $\tau$ as in type $C$ for a vexillary (signed) permutation $w=w(\tau)\in W_n$. The {\it covexillary Schubert variety} $X_{w_\circ w}$ in $Fl^B(V)$, the partition $\lambda$ associated to a covexillary (signed) permutation $w_\circ w$, $\mu$ associated to the pair $(w_\circ w,w_\circ v)$ with any $v\geq w$ are obtained in the same way as type $C$.

\section{Covexillary Formulas for Kazhdan-Lusztig Polynomials}\label{sec44}
\subsection{Relations to Grassmanianns of classical types}\label{gss}
In this section we review the relation of (co)vexillary elements with some inverse grassmannian ones in \cite{AIJK}. Here, an inverse grassmannian element is an element that is the inverse of a Grassmannian element in the Weyl group.

 Let ${P}$ be a maximal parabolic subgroup which contains the Borel subgroup ${B}$ of ${G}$ of rank $n$. Then ${G}/{P}$ has a structure of a homogeneous projective space. Let ${W}_{P}$ be the Weyl group of ${P}$ in the Weyl group ${W:=W_n}$ of $G$. The ${T}$-fixed points $({G}/{P})^{{T}}$ in ${G}/{P}$ are indexed by
the set of minimal representatives ${W}^{{P}}$ of the coset space ${W}/{W}_{{P}}$. For the classical Lie groups $G$, each element $w$ of ${W}^{{P}}$ can be represented by certain partitions. Specifically it can be represented by a weakly decreasing sequence $\rho=(\rho_1\geq\ldots\geq \rho_r> 0)$ for some $r$ for type $A$, and strictly decreasing sequence $\rho=(\rho_1>\cdots>\rho_r>0)$ for types $B$, $C$, and $D$. We express such elements as $w_\rho\in {W}^{{P}}$. One may refer \cite{BB} for these standard notations.

Let ${\overline{G}:}$ $=G(2n)$ be a simple linear algebraic group of classical type of rank $2n$ over $\mathbb{C}$. For the remaining of this article we use the notation $G$ for $G(n)$ of rank $n$ and $\overline{G}$ for $G(2n)$ of rank $2n$ to avoid any confusion that might occur, since we will deal with two different flag varieties. 

 Let $w_\nu\in\overline{W}^{\overline{P}}$ be a grassmannian element represented by the partition $\nu$ defined below for each classical type. Let $\lambda$ be the partition associated to a covexillary (signed) permutation $w_\circ w(\tau)$ with a triple $\tau=(\bf{k},\bf{p},\bf{q})$. 
 \begin{enumerate}
 \item Type $A$: Let $\overline{G}=SL(2n)$. Let $\overline{P}=P_{k_s}$ be the maximal parabolic subgroup associated to the $d$-th simple root where the simple roots are indexed as in \cite[p. 207]{BL}. 
  The grassmannian element $w_\nu\in \overline{W}^{\overline{P}}$ is represented by the partition $\nu=(\nu_1\geq \nu_2\geq\cdots\geq \nu_{k_s})$ defined by \eqref{nuA}.

 \item Type $C$: Let $\overline{G}$$=Sp(4n)$. The maximal parabolic subgroup $\overline{P}$$=P_{2n}$ of $\overline{G}$ is the one corresponding to the omission of the right end root, $2n$-th simple root in \cite[p. 207]{BL}. The partition $\nu$ is strictly decreasing defined by $\nu_{k_i}=3n+1-\lambda_{k_i}$ for all $i$ and $\nu_k=\nu_{k_i}+k_i-k$ for $k_{i-1}<k\leq k_{i}$.
 
 \item Type $D$: Here $\overline{G}=SO(4n)$ and $\overline{P}=P_{2n}$ is the maximal parabolic subgroup associated to the $2n$-th simple root as before. The partition $\nu$ is given by $\nu_{k_i}=3n-\lambda_{k_i}$ for $1\leq i\leq s$ and $\nu_{k}=\nu_{k_i}+k_i-k$ for $k_{i-1}< k\leq k_i$.
  
 \item Type $B$: The algebraic group is $\overline{G}=SO(4n+1)$. Let $\overline{P}=P_{2n}$ be the maximal parabolic subgroup that omits the $2n$-th simple root as in type $C$. Also, the partition $\nu$ is defined in the same way as in type $C$.
 \end{enumerate}

Let $w_\xi\in \overline{W}^{\overline{P}}$ be the element represented by a partition $\xi$ where $\xi$ is defined similarly as above with $\mu$ and the weak triple $\tau'=(\bf{k}',\bf{p},\bf{q})$ in \S\ref{sec4.1}.

We define the opposite Schubert cell by $\Omega_{w}^\circ=B^{-}wB/B$ and the opposite Schubert variety to be the closure of $\Omega_{w}^\circ$ in a flag variety $G/B$ for $w$ in the Weyl group. Theorem \ref{Thy} is a heart of the proof for Theorem \ref{M2}. The following theorem is \cite[Theorem 7.1]{AIJK} by specializing the Weyl group elements to our situation.

\begin{theorem}\label{Thy}
Let $w:=w({\tau})$ and $v:=w(\tau')$ be vexillary (signed) permutations in $W$, $w\leq v$. Let $\nu$ and $\xi$ be as above. For $w_\nu^{-1}$ and $v_\xi^{-1}$ in $\overline{W}$ with $w_\nu^{-1}\leq v_\xi^{-1}$, we have a scheme-theoretic isomorphism
\begin{equation}\label{isomorphism}
G/B=Fl^X(V)\supset\Omega_w\cap X_{v}^\circ\cong \Omega_{w_\nu^{-1}}\cap X_{{v_\xi}^{-1}}^\circ\subset Fl^X(V\oplus V)=\overline{G}/\overline{B}.
\end{equation}
Here $X=A,B,C$ or $D$. 
\end{theorem}

We note that the isomorphism \eqref{isomorphism} holds induced from the direct sum embedding 
\[
\sum: Fl^X(V)\rightarrow Fl^X(V\oplus V)
\]
in \cite[\S 6]{AIJK}. In particular, the image of $\Omega_w\cap X_v^\circ$ under the direct sum embedding $\sum$ is given by $\sum(\Omega_w\cap X_v^\circ)\cong \sum(X_v^\circ)\cap \Omega_{w_\nu^{-1}}$, and then certain matrix manipulations of $\sum(X_v^\circ)$ in \cite[\S 7]{AIJK} establishes $\sum(X_v^\circ)\cong X_{w_\xi^{-1}}^\circ$.

 A {\it Kazhdan-Lusztig variety}, so named in \cite{WY08}, is the intersection of a Schubert variety and the opposite Schubert cell in a flag variety. 
 For any $w\in W$ the opposite Schubert variety
$\Omega_w$ is isomorphic to the Schubert variety $X_{w_\circ w}$ \cite[\S 1.2]{KL}. We use $w_\circ$ for the longest element in both $W$ and $\overline{W}$ for abuse of notation, but it should be read in different Weyl groups. 
We will consider the Kazhdan-Lusztig varieties of the form $X_{ w}\cap \Omega_{ v}^\circ$ for $v\leq w$ for any $w,v\in W$. 
Then the isomorphism (\refeq{isomorphism}) can be interpreted as an isomorphsim
\begin{equation}\label{Equiv}
 X_{w_\circ w}\cap \Omega_{w_\circ v}^\circ\cong X_{w_\circ w_{\nu}^{-1}}\cap \Omega_{w_\circ v_{\xi}^{-1}}^\circ
\end{equation}
of Kazhdan-Lusztig varieties
for vexillaries $w, v\in W$ with $v\leq w$ and inverse grassmannians $w_\nu^{-1}, v_\xi^{-1}\in \overline{W}$ with $w_\nu^{-1}\leq v_\xi^{-1}$. 
Accordingly the isomorphism (\refeq{Equiv}) is an isomorphism of Kazhdan-Lusztig varieties in different flag varieties. It is noteworthy that the Kazhdan-Lusztig variety $X_{w_\circ w}\cap \Omega_{w_\circ v}^\circ$ associated to vexillary elements $w,v$ is in $Fl^X(V)=G/B$ while $X_{w_\circ w_{\nu}^{-1}}\cap \Omega_{w_\circ v_{\xi}^{-1}}^\circ$ associated to inverse grassmannians $w_\nu^{-1}, v_\xi^{-1}$ lies inside $Fl^X(V\oplus V)=\overline{G}/\overline{B}$.

\subsection{Isomorphisms and Main theorem } In this section, we show Theorem \ref{M2}.
Let us denote by $\mathbb{H}^*(Y, \mathcal{IC}(Y))$ the {\it hypercohomology group} of Deligne complex $\mathcal{IC}({Y})$ in $D^b(Y)$. The following proposition provides a relation between the hypercohomology and the stalk cohomology sheaf for the Kazhdan-Lusztig variety. It will allow us to relate \eqref{isomorphism} with the coefficients from \eqref{Poincare}, and eventually get  Theorem \ref{M2}. 

According to \cite[Lemma 4.5]{KL}, if $Z$ is a closed subvariety of $\mathbb{C}^N$ for some $N\geq 2$ invariant under the $(\mathbb{A}^1_{\mathbb{C}})^{\times}$-action by $t\cdot(z_1,\ldots,z_s)=(t^{\alpha_1}z_1,\ldots, t^{\alpha_s}z_s)$ for $\alpha_i>0$, $i=1,\ldots, s$, then
$\mathbb{H}^i(Z)=\mathcal{H}_{\mathbf{z}}^i(Z)$ for all $i$ where $\mathbf{z}$ is the cone point of $\mathbb{C}^N$. Furthermore, the Kazhdan-Lusztig varieties $X_{x}\cap \Omega^\circ_{y}$ be  for $x$ and $y$ in $W_n$ can be considered as $(\mathbb{A}^1_{\mathbb{C}})^{\times}$-stable closed subvarieties of $\mathbb{C}^{\ell(w_\circ)-\ell(y)}$ \cite[\S 1.5]{KL}. 
To be specific, $X_{x}\cap \Omega^\circ_{y}$ is a $T$-stable {\it attractive slice} to $X^\circ_{y}$ at $y$ in $X_x$. See \cite[\S 2.1]{Brion} for the definition of the attractive slice. That is, $y$ is an attractive fixed point for the $T$-action on $X_{x}\cap \Omega^\circ_{y}$. So, by \cite[Proposition A2]{Brion}, there exists a one parameter subgroup $\alpha: (\mathbb{A}^1_{\mathbb{C}})^{\times}\rightarrow T$ such that $\lim_{t\rightarrow0} \alpha(t)\cdot  z = y$ for all $z$ in a neighborhood of $y$. Moreover, the set of all $z\in X_{x}\cap \Omega^\circ_{y}$ satisfying $\lim_{t\rightarrow0} \alpha(t)\cdot  z = y$ forms an open affine $T$-stable neighborhood of $y$. Since $y\Omega_{id}$ is an open neighborhood of $y$ in $Fl$ and we restrict it on $X_{x}\cap \Omega^\circ_{y}$, all $z\in X_{x}\cap \Omega^\circ_{y}$ satisfies $\lim_{t\rightarrow0} \alpha(t)\cdot z = y$.
Note that $\chi_i \circ \alpha:(\mathbb{A}^1_{\mathbb{C}})^{\times}\rightarrow (\mathbb{A}^1_{\mathbb{C}})^{\times}$ is given by $t\rightarrow t^{\alpha_i}$ where $\chi_i$ are characters of $T$. In addition, by \cite[\S 1.3]{Brion} $Y:=X_{x}\cap \Omega^\circ_{y}$ can be embedded into $T_yY$ and all the weights of the $(\mathbb{A}^1_{\mathbb{C}})^{\times}$-action on $T_yY$ are positive, i.e., $\alpha_i>0$ for all $i$. Thus we have the following proposition.

\begin{proposition} \label{stalk}
Given $X_{x}\cap \Omega^\circ_{y}$ for $x$ and $y$ in $W$ with $x\geq y$ in Bruhat order, we have
 $$\mathbb{H}^i(X_{x}\cap \Omega^\circ_{y},\mathcal{IC}(X_{x}\cap \Omega^\circ_{y}))\cong\mathcal{H}_{ y}^i(\mathcal{IC}(X_{x}\cap \Omega^\circ_{y}))$$
  for all $i$.
 \end{proposition}

Now, we move toward Lemma \ref{basechange}. This lemma is an application of Theorem \ref{Thy}. Let us denote by
$i:y\Omega_{id}^\circ\hookrightarrow Fl$ the open immersion, where $y\Omega^\circ_{id}$ is the open neighborhood of a point $p_y$. 
Let $X_x=X_\ell\supset X_{\ell-2}\supset\cdots\supset X_0\supset \emptyset$ be the stratification of $X_x$ by closed subschemes of the Schubert variety $X_x$ of dimension $\ell:=\ell(x)$. Then the closed embedding $j:X_x\hookrightarrow Fl$ induces a stratification on $Fl$. So, with this stratification and \eqref{def:IC}, we obtain the Deligne sheaves $\mathcal{IC}({X_x}):=j_{!}\mathcal{IC}({X_x})$ on $Fl$. Moreover, with the similar argument and the map $i$, we have $\mathcal{IC}({X_x\cap y\Omega_{id}^\circ})$ from $i^{-1}\mathcal{IC}({X_x})\cong\mathcal{IC}({i^{-1}X_x})$.

The following lemma shows that $f: X_{x}\cap\Omega_{y}^\circ\rightarrow X_{x'}\cap\Omega_{y'}^\circ$ induces an isomorphism of hypercohomologies of $\mathcal{IC}$-sheaves on Kazhdan-Lusztig varieties.

\begin{lemma}\label{basechange}
Suppose that $f: X_{x}\cap\Omega_{y}^\circ\rightarrow X_{x'}\cap\Omega_{y'}^\circ$ is a scheme-theoretic isomorphism for pairs $x,y\in W$ and $x',y'\in \overline{W}$ with $y\leq x$ and $y'\leq x'$. We have 
$$\mathbb{H}^*(X_x\cap\Omega_y^\circ, \mathcal{IC}({X_{x}\cap\Omega_{y}^\circ}))\cong\mathbb{H}^*(X_{x'}\cap\Omega_{y'}^\circ,\mathcal{IC}({X_{x'}\cap\Omega_{y'}^\circ})).$$ 
\end{lemma}
\begin{proof}

Let $y\Omega_{id}^\circ=X_y^\circ\times\Omega_{y}^\circ$ be the open neighborhood of a point $p_y\in X_x$. As discussed earlier, we can get a $\mathcal{IC}$-sheaf on the open neighborhood ${X_x\cap y\Omega_{id}^\circ}$ of $p_y$ in $X_x$. 
In addition, for $x,y\in W$ with $x\leq y$, we know that $X_x\cap y\Omega_{id}^\circ\cong X_y^\circ\times (X_x\cap \Omega_y^\circ)$ by \cite[\S1.4]{KL}.

We define a $T$-equivariant isomorphism
 \begin{equation}\label{eqn:T-equi}
 \phi:=id^{-1}\times f^{-1}:X_{y}^\circ\times(X_{x'}\cap\Omega_{y'}^\circ)\rightarrow X_{y}^\circ\times(X_{x}\cap\Omega_{y}^\circ).
 \end{equation}
 
Then by topological invariance of $\mathcal{IC}$-sheaves \cite[p.97]{B}, we have
$$\phi^*\mathcal{IC}({X_{y}^\circ\times(X_{x}\cap\Omega_{y}^\circ)})\cong\mathcal{IC}({\phi^{-1}(X_y^\circ\times(X_{x}\cap\Omega_y^\circ))}).
$$
Thus, by \cite[Theorem 4.3.1]{B}, we obtain $\mathcal{IC}({\phi^{-1}(X_y^\circ\times(X_{x}\cap\Omega_y^\circ))}) \cong\mathcal{IC}({X_{y}^\circ\times(X_{x'}\cap\Omega_{y'}^\circ)})$. As $\phi$ is an isomorphism, we get
 $$\mathcal{IC}({X_{y}^\circ\times(X_{x}\cap\Omega_{y}^\circ)})\cong\mathcal{IC}({X_{y}^\circ\times(X_{x'}\cap\Omega_{y'}^\circ)}).$$
Then we have
  $$\mathcal{IC}({X_y^\circ})\boxtimes\mathcal{IC}({X_x\cap\Omega_y^\circ})\cong\mathcal{IC}({X_{y}^\circ})\boxtimes\mathcal{IC}({X_{x'}\cap\Omega_{y'}^\circ})$$
 from \cite[\S 5]{L01} where $\boxtimes$ is the geometrical external tensor product of perverse sheaves.
 
As $X_y^\circ$ is smooth and $\mathcal{IC}$ sheaf is constructible, its intersection cohomology sheaf is constant on $X_y^\circ$. Hence, we conclude $\mathcal{IC}({X_x\cap \Omega_y^\circ})\cong\mathcal{IC}({X_{x'}\cap\Omega_{y'}^\circ})$ of Deligne complexes of sheaves, which immediately gives the statement.
\end{proof}

\begin{remark}
It is expected that if $f:A\rightarrow B$ is a scheme-theoretic isomorphism, then $\mathcal{IC}(A)\cong \mathcal{IC}(B)$. However, to the best of the author's knowledge, this is only true if there exists a $\mathcal{IC}$-sheaf on $B$, and $A$ and $B$ admit stratifications which are homeomorphic, according to \cite[p.97]{B}. To resolve this, in Lemma \ref{basechange}, we pull back the stratification of $X_x\cap y\Omega_{id}^\circ$ under the isomorphism \eqref{eqn:T-equi}, and then get a $\mathcal{IC}$-sheaf on the Kazhdan-Lusztig varieties ${X_{x'}\cap\Omega_{y'}^\circ}$ and ${X_x\cap\Omega_y^\circ}$ as in \cite[Lemma 3.2]{P}.

\end{remark}

The following proposition explains that  two relevant $T$-fixed points are included in a Schubert variety. We consider the locally closed set $\mathcal{Z}_{\tau'}=\{E_\bullet\;|\;\mathrm{dim}(E_{p_i}\cap F_{q_i})=k_i'\}\subset Fl=G/B$ where $\tau'=(\bf{k}',\bf{p},\bf{q})$ is a weak triple with a vexillary element $v'=w(\tau')$ in the Weyl group $W$. Suppose that the stabilizer of the action on the partial flag $F_{q_i}$ is $Q$. Then the set $\mathcal{Z}_{\tau'}$ is a $Q$-orbit of $p_{w_\circ v'}$ as a disjoint union of Schubert cells containing both $T$-fixed points $p_{w_\circ v}$ and $p_{w_\circ v'}$. That is,
\begin{proposition}
Let $w$ be a vexillary element. A Schubert variety $X_{w_\circ w}$ contains both $p_{w_\circ v}$ and $p_{w_\circ v'}$ for any $v$, and a vexillary (signed) permutation $v'(\tau')$ induced by the weak triple in $W$.  
\end{proposition}
\begin{proof}
We fix $g\in G$ satisfying $gp_{w_\circ v}=p_{w_\circ v'}$. 
Essentially $g$ stabilizes the partial flag $F_{q_1}\subseteq \cdots \subseteq F_{q_s}$, which implies that $\mathcal{Z}_{\tau'}$ is preserved by the multiplication by $g$, that is, $g\mathcal{Z}_{\tau'}=\mathcal{Z}_{\tau'}$. For a similar reason, $X_{w_\circ w}$ is invariant under the action of $g$ so that $gX_{w_\circ w}=X_{w_\circ w}$. One can check by the invariance under $g$-action that $g(X_{w_\circ w}\cap\mathcal{Z}_{\tau'})=gX_{w_\circ w}\cap g\mathcal{Z}_{\tau'}=X_{w_\circ w}\cap \mathcal{Z}_{\tau'}$. This completes the proposition.
\end{proof}
We observe that $p_{w_\circ v}$ and $p_{w_\circ v'}$ can be in different strata of $X_{w_\circ w}$.

Let $X$ and $Y$ be varieties with points $x\in X$ and $y\in Y$. We say that the singularity of $X$ at $x$ and the singularity of $Y$ at $y$ are {\it smoothly equivalent} if there is a variety $Z$ such that two maps $f_1:Z\rightarrow X$ and $f_2:Z\rightarrow Y$ 
are smooth at $z$ with $f_1(z)=x$ and $f_2(z)=y$, \cite[Definition 3.7]{Juteau}.

Finally, we are in the position to show Theorem \ref{MainTh}. 
\begin{theorem}
\label{MainTh}
For a covexillary Schubert variety $X_{w_\circ w}$ and any torus-fixed point $p_{w_\circ v}\in X_{w_\circ w}$ in $Fl=G/B$, there exist $X_{w_\circ w_\nu^{-1}}$ and a torus-fixed point $p_{w_\circ v_\xi^{-1}}\in X_{w_\circ w_\nu^{-1}}$ in $\overline{G}/\overline{B}$ such that
\[
\mathcal{H}_{p_{w_\circ v}}^{2q}(\mathcal{IC}(X_{w_\circ {w}}))=\mathcal{H}_{p_{w_\circ v_\xi^{-1}}}^{2q}(\mathcal{IC}(X_{w_\circ w_\nu^{-1}}))
\]
for all $q\geq 0$.
In particular, $P_{w_\circ {v},w_\circ w}(q)=P_{w_\circ {v}_\xi^{-1}, w_\circ w_\nu^{-1}}(q).$
\end{theorem}

\begin{proof}

Let $m_x$ be the multiplication map by $x\in G$ and $e$ the identity of $G$.
We take the commutative diagram
$$
\begin{tikzcd}
X_{w_\circ w} &\arrow{l}{}[swap]{f_1} X_{w_\circ w}\cap p_{w_\circ v'}\Omega_{id}^\circ&\arrow{l}{}[swap]{m_g} X_{w_\circ w}\cap p_{w_\circ v}\Omega_{id}^\circ\arrow[bend right=20]{ll}[name=O]{h} \arrow[r,"f_2"]& X_{w_\circ w}\\
p_{w_\circ v'}\arrow[u, "u_1"]&\arrow[l,"m_e"]p_{w_\circ v'}\arrow[u, "u_1"]&\arrow[l,"m_g"]\arrow[bend left=20]{ll}[swap]{m_g} p_{w_\circ v}\arrow{r}{}[swap]{m_e} \arrow[u, "u_2"]& p_{w_\circ v}\arrow[u, "u_2"]
\end{tikzcd},
$$
where $u_i$ is the inclusion, $f_i$ is an open immersion, and $h=f_1\circ m_g$ is the composition of $f_1$ and $m_g$. Since $f_i$ is smooth at $p_{w_\circ v}$, the singularity of $X_{w_\circ w}$ at $p_{w_\circ v'}$ is smoothly equivalent to the singularity of $X_{w_\circ w}$ at $p_{w_\circ v}$. By virtue of \cite[Proposition 3.8]{Juteau}, we deduce 
$\mathcal{IC}(X_{w_\circ w})_{p_{w_\circ v}}\cong \mathcal{IC}(X_{w_\circ w})_{p_{w_\circ v'}}$
so that
\begin{equation}\label{isostalk}
\mathcal{H}_{p_{w_\circ v}}^{2q}(\mathcal{IC}(X_{w_\circ {w}}))\cong\mathcal{H}_{p_{w_\circ v'}}^{2q}(\mathcal{IC}(X_{w_\circ w})).
\end{equation}
The definition of stalk for the intersection cohomology sheaf shows the isomorphism 
\begin{equation}\label{00}
\mathcal{H}_{p_{w_\circ v'}}^{2q}(\mathcal{IC}(X_{w_\circ w}))\cong \mathcal{H}_{p_{w_\circ v'}}^{2q}(\mathcal{IC}(X_{w_\circ w}\cap \Omega^\circ_{{w_\circ v'}})).
\end{equation}
We know that $w=w(\tau)$ and $v'$ are vexillary elements in $W$. Combined with (\refeq{isostalk}), Proposition \ref{stalk}, (\refeq{Equiv}) and Lemma \ref{basechange}, we obtain 
\begin{align*}
\mathcal{H}_{p_{w_\circ v}}^{2q}(\mathcal{IC}(X_{w_\circ {w}}))&\cong\mathbb{H}^{2q}(X_{w_\circ w_\nu^{-1}}\cap \Omega_{w_\circ v_\xi^{-1}}^\circ, \mathcal{IC}(X_{w_\circ w_\nu^{-1}}\cap \Omega_{w_\circ v_\xi^{-1}}^\circ))\\
&\cong \mathcal{H}_{p_{w_\circ v_\xi^{-1}}}^{2q}(\mathcal{IC}(X_{w_\circ w_\nu^{-1}}\cap \Omega^\circ_{{w_\circ v_\xi^{-1}}})).
\end{align*}
Since $\mathcal{H}_{p_{w_\circ v_\xi^{-1}}}^{2q}(\mathcal{IC}(X_{w_\circ w_\nu^{-1}}))\cong  \mathcal{H}_{p_{w_\circ v_\xi^{-1}}}^{2q}(\mathcal{IC}(X_{w_\circ w_\nu^{-1}}\cap \Omega^\circ_{{w_\circ v_\xi^{-1}}}))$, we have thereby proved that their Hilbert-Poincar\'{e} polynomials coincide. 
 
\end{proof}

Finally, we finish a proof of our main results by the following theorem.
 \begin{theorem}\label{M2}
Let ${w_\circ w}$ be a covexillary element in $W$, and $v$ be any element with $w\leq v$. Then there are some grassmannian elements $w_\nu, v_\xi\in \overline{W}^{\overline{P}}$ with $w_\nu\geq v_\xi$ such that
$$P_{w_\circ {v},w_\circ w}(q)=P_{w_\circ {v}_\xi, w_\circ w_\nu}(q).$$
\end{theorem}
   \begin{proof}
We know from \cite[Cor. 4.4]{Brenti}, \cite[Prop. 7.6 and Sec. 7.13]{Hum} that $P_{y,w}(q)=P_{y^{-1},v^{-1}}(q)=P_{w_\circ y w_\circ,w_\circ w w_\circ}(q)$ for $u,v\in W$. 
Since 
$$P_{w_\circ {v}_\xi^{-1}, w_\circ w_\nu^{-1}}(q)=P_{{v}_\xi w_\circ^{-1},w_\nu w_\circ^{-1}}(q)=P_{w_\circ {v}_\xi,w_\circ w_\nu}(q),$$
we conclude the statement with Theorem \ref{M2}.
\end{proof}

\subsection{Combinatorial formulas for covexillary cases}\label{sec1}
In this section, we discuss the combinatorial and inductive formulas for Kazhdan-Lusztig polynomials in covexillary case, generalizing the cograssmannian cases in classical types.

Let $w_\circ w(\tau)$ be a covexillary (signed) permutation with a triple $\tau$, and $w_\circ v$ be any element with $w(\tau)\leq v$. Recall from \S\ref{sec4.1} that $\lambda(w_\circ w(\tau))$ and $\mu(w_\circ w(\tau), w_\circ v)$ be the $\alpha\beta$-words.  

Given the $\alpha\beta$-words $\lambda(w_\circ w(\tau))$ and $\mu(w_\circ w(\tau), w_\circ v)$, one can define the set 
\begin{equation}\label{tree}
A^{G}(\lambda(w_\circ w(\tau))/\mu(w_\circ v))
\end{equation}
 of all labelings of tree(s) $\mathcal{T}:=\mathcal{T}(\lambda(w_\circ w(\tau)))$ for $G=SL(n)$ of type $A$ by \cite[6.5]{LS} and for $G=Sp(2n)$ of type $C$, $G=SO(2n)$ of type $D$ and $G=SO(2n+1)$ of type $B$ by \cite{B88}. We note that the construction of tree(s) $\mathcal{T}$ and the set $A^G$ of labelled tree(s) in \cite{LS,B88} arising only from $\alpha\beta$-words. So, the set \eqref{tree} and tree(s) $\mathcal{T}$ are well-defined with $\alpha\beta$-words.

 The following gives the combinatorial formulas for Kazhdan-Lusztig polynomials associated to covexillary elements.  
 
 \begin{corollary}\label{MainI}
Let $G$ be an algebraic group of classical types. Let $w:= w(\tau)$ be a vexillary element with a triple $\tau$ in $W$. Let $v$ be any element in $W$ with $w\leq v$ in Bruhat order.  Then the Kazhdan-Lusztig polynomial is given by 
$$P_{w_\circ v, w_\circ w}(q)=\sum q^{|\mathcal{T}|}$$
where the sums take over all edge labelings of $\mathcal{T}$ in $A^G(\lambda(w_\circ w),\mu(w_\circ w,w_\circ v))$. $|\mathcal{T}|$ denotes the sum of all the edge labels.
\end{corollary}
\begin{proof}
Let $w_\nu$ and $w_\xi$ be the grassmannian elements defined in \S\ref{gss}. Since the corresponding $\alpha\beta$ words of cograssmannian elements $w_\circ w_\nu$ and $w_\circ {v}_\xi$ are $\lambda(w_\circ w)$ and $\mu(w_\circ w, w_\circ v)$ respectively, by \cite[Theorem 7.8]{LS} for type $A$ and \cite{B88} for types $C$, $D$ and $B$, we have
\[
P_{w_\circ {v}_\xi,w_\circ w_\nu}(q)=\sum q^{|\mathcal{T}|}.
\]
Combined with Theorem \ref{M2}, we establish the result.
\end{proof}

We remark that the Kazhdan-Lusztig polynomials depend only on triples $\tau$ and $\tau'$ from vexillary element $w$ and any $v$ with $w\leq v$ in $W$. In particular if $v'\geq w$ is another element with $\tau'(v',w)=\tau(v,w)$, the Kazhdan-Lusztig polynomials are equal. In other words, the $\alpha\beta$-words for covexillary cases arise from the sequences of $\tau$ and $\tau'$, while the ones for cograssmannian cases are only from the partitions.

\subsection{\bf{Covexillary version of Inductive formulas for Kazhdan-Lusztig polynomials}}\label{Main result II}

Recall that $\lambda$ is the partition associated to a covexillary (signed) permutation $w_\circ w(\tau)$ with a triple $\tau=(\bf{k},\bf{p},\bf{q})$ and $\mu$ be the partition associated to the pair $(w_\circ w, w_\circ v)$ for any element $v$ with $w\leq v$.

The algorithm of this second result can be realized by two steps with the matrix
$
\mathfrak{H}(\lambda)=\begin{pmatrix}
a_1 & \cdots& a_d\\
b_0 &\cdots& b_{d-1}\\
\end{pmatrix}_\lambda
$ defined in \S\ref{sec4.1}. First, we find the smallest $i$ such that $b_i\leq a_i$ and $a_{i+1}\leq b_{i+1}$. Next we replace 
$\label{ind}
\mathfrak{H}(\lambda)=\begin{pmatrix}
\cdots &a_i& a_{i+1}&a_{i+2}&\cdots\\
\cdots &b_{i-1}&b_i& b_{i+1}&\cdots\\
\end{pmatrix}_\lambda
$ by
$$
\mathfrak{H}_1:=\begin{pmatrix}
\cdots &a_i+ a_{i+1}&a_{i+2}&\cdots\\
\cdots &b_{i-1}&b_i+ b_{i+1}&\cdots\\
\end{pmatrix}.
$$
Note that $d=s$ or $d=s+1$, depending on their types, i.e., type $D$ could have $d=s+1$ for a triple $\tau$ of length $s$.

In addition we impose two conditions on $\mathfrak{H}(\lambda)$ by 
\begin{enumerate}\label{twoconditions}
\item $\lambda_d<N-a_d$
\item ($a_d+a_{d-1}+\cdots+a_i)-(b_{d-1}+\cdots+b_i)<N-\lambda_d,\quad\text{for}\quad 1\leq i$,
\end{enumerate}
where $\lambda_d=p_d+q_d$ and $N=2n+1$ if $G=Sp(2n)$ or $SO(2n+1)$, and $\lambda_d=p_d+q_d+2$ and $N=2n$ if $G=SO(2n)$. There are no such hypotheses on type $A$. 

Given $\mu\leq \lambda$, $c(\mu,\lambda)$ is the {\it depth} sequence of $\mu$ in $\lambda$ such that
\[
\mathfrak{H}(\mu)=\begin{pmatrix}
e_1 &e_2& \cdots& e_m&e_{m+1}\\
f_0 &f_1& \cdots& f_{m-1}&f_m\\
\end{pmatrix}_\mu
=\begin{pmatrix}
a_1+c_1 &a_2+c_2-c_1& \cdots& a_m+c_m-c_{m-1}&1\\
f_0 &f_1& \cdots& f_{m-1}&f_m\\
\end{pmatrix}_\mu
\]
where the last column appears precisely when $G=SO(2n)$ and $c_m$ is odd, and $2n=\sum_{i=1}^{m+1} e_i+\sum_{i=0}^{m}f_i$ \cite[3.2]{SV}.

Let
$$\genfrac{[}{]}{0pt}{}{\alpha}{\beta}=\dfrac{(q^\alpha-1)\cdots(q^{\alpha-\beta+1}-1)}{(q^\beta-1)(q^{\beta-1}-1)\cdots(q-1)}$$
 where $\genfrac{[}{]}{0pt}{}{\alpha}{\beta}=0$ if $\beta<0$ or $\alpha<\beta$, and $\genfrac{[}{]}{0pt}{}{0}{0}=1$ as well as $\genfrac{[}{]}{0pt}{}{\alpha}{\beta}=\genfrac{[}{]}{0pt}{}{\alpha}{\alpha-\beta}$.

Our inductive formula can be presented as follows.

\begin{corollary}\label{MainII}
Let $w=w(\tau)$ be a vexillary (signed) element in $W$. Let $v$ be any element with $w\leq v$. The inductive formula for Kazhdan-Lusztig polynomial associated to $w_\circ w$ and $w_\circ v$ is given by
$$P_{w_\circ v,w_\circ w}(q)=\sum_tq^{(c_i-t)(c_{i+1}-t)}\genfrac{[}{]}{0pt}{}{a_{i+1}-c_i+c_{i+1}}{c_{i+1}-t}\genfrac{[}{]}{0pt}{}{b_i+c_i-c_{i+1}}{c_{i}-t}P_{\mathfrak{H}_1,\bold{c}(t)}(q),$$
where $\bold{c}(t)=(c_0,\ldots,c_{i-1},t,c_{i+2},\ldots, c_d)$ for a capacity $\bold{c}(\mu,\lambda)=(c_0, \ldots, c_d)$ and $c_0=0$.
\end{corollary}
\begin{proof}
Let $w_\nu$ and $w_\xi$ be the grassmannian elements defined in \S\ref{gss}. Since the corresponding $2\times m$ matrix of cograssmannian elements $w_\circ w_\nu$ and $w_\circ {v}_\xi$ are $\mathfrak{H}(\lambda(w_\circ w))$ and $\mathfrak{H}(\mu(w_\circ w, w_\circ v))$ respectively, by \cite[Theorem 2]{Z} for type $A$ and \cite{SV} for types $C$, $D$ and $B$ or \cite[Chapter 9.1.20]{BL}, we have
\[
P_{w_\circ {v}_\xi,w_\circ w_\nu}(q)=\sum_tq^{(c_i-t)(c_{i+1}-t)}\genfrac{[}{]}{0pt}{}{a_{i+1}-c_i+c_{i+1}}{c_{i+1}-t}\genfrac{[}{]}{0pt}{}{b_i+c_i-c_{i+1}}{c_{i}-t}P_{\mathfrak{H}_1,\bold{c}(t)}(q),
\]
where $\bold{c}(t)=(c_0,\ldots,c_{i-1},t,c_{i+2},\ldots, c_d)$ for a capacity $\bold{c}=(c_0, \ldots, c_d)$ and $c_0=0$.
Hence, the result follows by Theorem \ref{M2}.
\end{proof}

\begin{acknowledgments}
This paper is initiated by Alexander Woo's question in his private conversation with David Anderson. We are grateful to his piercing insight for suggesting this problem.
We are sincerely grateful to Daniel Juteau and Nicholas Proudfoot for kind explanations and references. We also wish to thank Alexander Yong for his keen interest in this topic, his multiple reviews of the manuscript, and his valuable feedback. 
Most importantly the author would like to express indebtedness to David Anderson for useful comments on the earlier draft and his generous guidance. Lastly, we thank the anonymous referees for their valuable comments, which have significantly improved this manuscript.
\end{acknowledgments}

\bibliographystyle{amsplain}

\begin{bibdiv}
\begin{biblist}



\bib{AF12}{article}{
   author={Anderson, David},
   author={Fulton, William},
   title={Degeneracy loci, Pfaffians, and vexillary signed permutations in types $B$, $C$,
and $D$},
   journal={preprint, arXiv:1210.2066},
   date={2012},
}

\bib{AF20}{article}{
   author={Anderson, David},
   author={Fulton, William},
   title={Vexillary signed permutations revisited},
   journal={Algebr. Comb.},
   volume={3},
   date={2020},
   number={5},
   pages={1041--1057},
}
		
	\bib{AIJK19}{article}{
   author={Anderson, Dave},
   author={Ikeda, Takeshi},
   author={Jeon, Minyoung},
   author={Kawago, Ryotaro},
   title={Multiplicities of Schubert varieties in the symplectic flag
   variety},
   language={English, with English and French summaries},
   journal={S\'{e}m. Lothar. Combin.},
   volume={82B},
   date={2020},
   pages={Art. 95, 12},
}		

\bib{AIJK}{article}{
   author={Anderson, David},
   author={Ikeda, Takeshi},
   author={Jeon, Minyoung},
   author={Kawago, Ryotaro},
   title={The multiplicity of a singularity in a vexillary Schubert variety},
   journal={Adv. Math.},
   volume={435},
   date={2023},
   pages={Paper No. 109366, 39},
}
%  



\bib{B}{book}{
   author={Banagl, M.},
   title={Topological invariants of stratified spaces},
   series={Springer Monographs in Mathematics},
   publisher={Springer, Berlin},
   date={2007},
   pages={xii+259}, 
}
		


\bib{BL}{book}{
   author={Billey, Sara},
   author={Lakshmibai, V.},
   title={Singular loci of Schubert varieties},
   series={Progress in Mathematics},
   volume={182},
   publisher={Birkh\"{a}user Boston, Inc., Boston, MA},
   date={2000},
   pages={xii+251},
  }
  \bib{BilleyLam}{article}{
   author={Billey, Sara},
   author={Lam, Tao Kai},
   title={Vexillary elements in the hyperoctahedral group},
   journal={J. Algebraic Combin.},
   volume={8},
   date={1998},
   number={2},
   pages={139--152},
}
\bib{BB}{book}{
   author={Bj\"{o}rner, Anders},
   author={Brenti, Francesco},
   title={Combinatorics of Coxeter groups},
   series={Graduate Texts in Mathematics},
   volume={231},
   publisher={Springer, New York},
   date={2005},
   pages={xiv+363},
}

\bib{B88}{article}{
   author={Boe, Brian D.},
   title={Kazhdan-Lusztig polynomials for Hermitian symmetric spaces},
   journal={Trans. Amer. Math. Soc.},
   volume={309},
   date={1988},
   number={1},
   pages={279--294},
}	



	
\bib{Brenti}{article}{
   author={Brenti, Francesco},
   title={A combinatorial formula for Kazhdan-Lusztig polynomials},
   journal={Invent. Math.},
   volume={118},
   date={1994},
   number={2},
   pages={371--394},
}
	
		
  \bib{Brion}{article}{
   author={Brion, M.},
   title={Rational smoothness and fixed points of torus actions},
   note={Dedicated to the memory of Claude Chevalley},
   journal={Transform. Groups},
   volume={4},
   date={1999},
   number={2-3},
   pages={127--156},
}
	
	




\bib{GM}{article}{
   author={Goresky, Mark},
   author={MacPherson, Robert},
   title={Intersection homology. II},
   journal={Invent. Math.},
   volume={72},
   date={1983},
   number={1},
   pages={77--129},
}
\bib{Hum}{book}{
   author={Humphreys, James E.},
   title={Reflection groups and Coxeter groups},
   series={Cambridge Studies in Advanced Mathematics},
   volume={29},
   publisher={Cambridge University Press, Cambridge},
   date={1990},
   pages={xii+204},
}





\bib{JW}{article}{
   author={Jones, Brant},
   author={Woo, Alexander},
   title={Mask formulas for cograssmannian Kazhdan-Lusztig polynomials},
   journal={Ann. Comb.},
   volume={17},
   date={2013},
   number={1},
   pages={151--203},
}



\bib{J1}{article}{
   author={Joshua, Roy},
   title={Vanishing of odd-dimensional intersection cohomology},
   journal={Math. Z.},
   volume={195},
   date={1987},
   number={2},
   pages={239--253},
}

\bib{Juteau}{article}{
   author={Juteau, Daniel},
   title={Decomposition numbers for perverse sheaves},
   language={English, with English and French summaries},
   journal={Ann. Inst. Fourier (Grenoble)},
   volume={59},
   date={2009},
   number={3},
   pages={1177--1229},
 }





\bib{KL0}{article}{
   author={Kazhdan, David},
   author={Lusztig, George},
   title={Representations of Coxeter groups and Hecke algebras},
   journal={Invent. Math.},
   volume={53},
   date={1979},
   number={2},
   pages={165--184},
}
\bib{KL}{article}{
   author={Kazhdan, David},
   author={Lusztig, George},
   title={Schubert varieties and Poincar\'{e} duality},
   conference={
      title={Geometry of the Laplace operator},
      address={Proc. Sympos. Pure Math., Univ. Hawaii, Honolulu, Hawaii},
      date={1979},
   },
   book={
      series={Proc. Sympos. Pure Math., XXXVI},
      publisher={Amer. Math. Soc., Providence, R.I.},
   },
   date={1980},
   pages={185--203},

}





	

  

	
	\bib{L}{article}{
   author={Lascoux, Alain},
   title={Polyn\^{o}mes de Kazhdan-Lusztig pour les vari\'{e}t\'{e}s de Schubert
   vexillaires},
   language={French, with English and French summaries},
   journal={C. R. Acad. Sci. Paris S\'{e}r. I Math.},
   volume={321},
   date={1995},
   number={6},
   pages={667--670},
}
\bib{LS}{article}{
   author={Lascoux, Alain},
   author={Sch\"{u}tzenberger, Marcel-Paul},
   title={Polyn\^{o}mes de Kazhdan \& Lusztig pour les grassmanniennes},
   language={French},
   conference={
      title={Young tableaux and Schur functors in algebra and geometry
      (Toru\'{n}, 1980)},
   },
   book={
      series={Ast\'{e}risque},
      volume={87},
      publisher={Soc. Math. France, Paris},
   },
   date={1981},
   pages={249--266},
}		
\bib{LY}{article}{
   author={Li, Li},
   author={Yong, Alexander},
   title={Kazhdan-Lusztig polynomials and drift configurations},
   journal={Algebra Number Theory},
   volume={5},
   date={2011},
   number={5},
   pages={595--626},
}


\bib{L01}{article}{
   author={Lyubashenko, V. V.},
   title={External tensor product of perverse sheaves},
   language={Ukrainian, with English and Ukrainian summaries},
   journal={Ukra\"{\i}n. Mat. Zh.},
   volume={53},
   date={2001},
   number={3},
   pages={311--322},
   translation={
      journal={Ukrainian Math. J.},
      volume={53},
      date={2001},
      number={3},
      pages={354--367},
    },
}
	

   \bib{M}{book}{
   author={L. Maxim},
   title={Intersection Homology and Perverse Sheaves with Applications to Singularities},
 note={Graduate Texts in Mathematics, No. 281},
  publisher={Springer-Verlag, New York-Heidelberg},
   date={2019},
   pages={xv+270},
}

	
\bib{P}{article}{
   author={Proudfoot, Nicholas},
   title={The algebraic geometry of Kazhdan-Lusztig-Stanley polynomials},
   journal={EMS Surv. Math. Sci.},
   volume={5},
   date={2018},
   number={1-2},
   pages={99--127},
}


\bib{SV}{article}{
   author={Sankaran, Parameswaran},
   author={Vanchinathan, P.},
   title={Small resolutions of Schubert varieties and Kazhdan-Lusztig
   polynomials},
   journal={Publ. Res. Inst. Math. Sci.},
   volume={31},
   date={1995},
   number={3},
   pages={465--480},

}

	
\bib{Tam}{article}{
   author={Tamvakis, Harry},
   title={Degeneracy locus formulas for amenable Weyl group elements},
   journal={preprint, arXiv:1909.06398},
   date={2019},
}
	


	\bib{WY08}{article}{
   author={Woo, Alexander},
   author={Yong, Alexander},
   title={Governing singularities of Schubert varieties},
   journal={J. Algebra},
   volume={320},
   date={2008},
   number={2},
   pages={495--520},
}
	\bib{WY23}{article}{
   author={Woo, Alexander},
   author={Yong, Alexander},
   title={Schubert geometry and combinatorics},
   journal={preprint, arXiv:2303.01436},
   date={2023},
}

	
\bib{Z}{article}{
   author={Zelevinsky, A. V.},
   title={Small resolutions of singularities of Schubert varieties},
   journal={Funktsional. Anal. i Prilozhen.},
   volume={17},
   date={1983},
   number={2},
   pages={75--77},
}
\end{biblist}
\end{bibdiv}

\end{document}